\documentclass[12pt,leqno]{article}
\usepackage{amsmath, amssymb}

\def\epsilon{\varepsilon}

\begin{document}

\LARGE \noindent 
{\bf A remark on embedding of a cylinder on a real 
commutative Banach algebra}

\large

\vspace*{0.4em}

\hfill Hiroki Yagisita (Kyoto Sangyo University)

\vspace*{1.6em}

\normalsize

Abstract: 

Let $A$ be a real commutative Banach algebra with unity.  
Let $a_0\in A\setminus\{0\}$. Let $\mathbb Z a_0:=\{na_0\}_{n\in \mathbb Z}$. 
Then, $\mathbb Z a_0$ is a discrete subgroup of $A$.  
For any $n\in \mathbb Z$, the Frechet derivative of the mapping  
$$x \, \in \, A \ \ \ \mapsto \ \ \ x+na_0 \, \in \, A$$
is the identity map on $A$ and, especially, an $A$-linear transformation on $A$. 
So, the quotient group $A/(\mathbb Z a_0)$ is a $1$-dimensional $A$-manifold 
and the covering projection 
$$x \, \in \, A \ \ \ \mapsto \ \ \ x+\mathbb Z a_0 \, \in \, A/(\mathbb Z a_0)$$ 
is an $A$-map. We call $A/(\mathbb Z a_0)$ 
the $1$-dimensional $A$-cylinder by $a_0$. 

Let $T$ be a compact Hausdorff space. Suppose that 
there exist $t_1\in T$ and $t_2\in T$ such that $t_1\not=t_2$ holds. 
Then, the set $C(T;\mathbb R)$ of all real-valued continuous functions on $T$ 
is a real commutative Banach algebra with unity 
and $\mathbb R \, \subsetneq \, C(T;\mathbb R)$ holds. 
In this paper, we show that there exists $a_0 \, \in \, C(T;\mathbb R)\setminus \mathbb R$ 
such that for any $k\, \in \, \mathbb N$, the $1$-dimensional $C(T;\mathbb R)$-cylinder 
$(C(T;\mathbb R))/(\mathbb Z a_0)$ by $a_0$
cannot be embedded in the finite direct product space $(C(T;\mathbb R))^k$ 
as a $C(T;\mathbb R)$-submanifold.

%Also, $f$ is said to be $\mathbb C^n$-triangular, 
%if $f$ is complex analytic and for any $i$ and $j$, 
%$i<j$ implies $\frac{\partial f_i}{\partial z_j}=0$. 
%Kasuya suggested that a $\mathbb C^n$-analytic manifold 
%and a $\mathbb C^n$-triangular manifold might, for example, 
%be related to a holomorphic web and a holomorphic foliation.  

%\vfill 
\vspace*{0.8em} 
\noindent 
%Keywords: \ 

%\noindent 
%immersion, bump function, partition of unity, infinite-dimensional Lie group. 


\vfill 
\vspace*{0.8em} 
\noindent 
The related literature: \ 
``https://www.researchgate.net/profile/Hiroki\_Yagisita''

%\vfill 
\noindent 
%Department of Mathematics, Faculty of Science, Kyoto Sangyo University, Motoyama, Kamigamo, Kita-ku, Kyoto-City, 603-8555, Japan. 

%\vspace*{0.4em} 

\noindent 
%e-mail: hrk0ygst@cc.kyoto-su.ac.jp

\newpage 

In this paper, we give a remark concerned with an $A$-manifold 
(a manifold on a commutative topological algebra $A$). 
There are already various studies related to $A$-manifolds 
or their analogues (e.g., [1], [2], $\cdots$, [9]). 
In [10], we showed the existence of a $C([0,1];\mathbb R)$-manifold 
that cannot be embedded 
in the finite-dimensional Affine space $(C([0,1];\mathbb R))^k$ 
as a $C([0,1];\mathbb R)$-submanifold. In this paper, 
we show the following, which is a generalization. 

\vspace*{0.8em} 

\noindent 
{\bf Theorem} : \ \ \ Let $T$ be a compact Hausdorff space. Suppose that 
there exist $t_1\in T$ and $t_2\in T$ such that $t_1\not=t_2$ holds. 
Then, there exists $a_0 \, \in \, C(T;\mathbb R)\setminus \mathbb R$ 
such that for any $k\, \in \, \mathbb N$, the cylinder $(C(T;\mathbb R))/(\mathbb Z a_0)$ 
cannot be embedded in the Cartesian product $(C(T;\mathbb R))^k$ 
as a $C(T;\mathbb R)$-submanifold. 

\vspace*{0.8em} 

\noindent 
{\sf Proof} : \ 
Because $T$ is normal, by Urysohn's lemma, there 
exists $a_0 \, \in \, C(T;\mathbb R)\setminus \mathbb R$ 
such that $a_0(t_1)=-1$ and $a_0(t_2)=+1$ hold. Let $k\in \mathbb N$. 
Then, by a contradiction, we show that $(C(T;\mathbb R))/(\mathbb Z a_0)$ 
cannot be embedded in $(C(T;\mathbb R))^k$ 
as a $C(T;\mathbb R)$-submanifold. 

Suppose that $(C(T;\mathbb R))/(\mathbb Z a_0)$ 
can be embedded in $(C(T;\mathbb R))^k$ 
as a $C(T;\mathbb R)$-submanifold. 
Then, there exists a $C(T;\mathbb R)$-injection $\Psi$ 
from $(C(T;\mathbb R))/(\mathbb Z a_0)$ to $(C(T;\mathbb R))^k$. Let 
\begin{equation}b_0(t) \, := \, \max \{0,a_0(t)\} \ \ \ \ \ \ (t\in T).\end{equation}
Then, $b_0\in C(T;\mathbb R)$ and 
$$[0]=[a_0]\not=[b_0]$$
hold, where $[x]$ denotes the equivalence class of $x \in C(T;\mathbb R)$. 
That is, we put  
$$[x] \ := \ x \, + \, \mathbb Z a_0 \ \ \ \ \ \ (x \in C(T;\mathbb R)).$$
So, because $\Psi([0])=\Psi([a_0])\not=\Psi([b_0])$ holds, 
there exists $l_0\in \{1,2,\cdots,k\}$ such that 
\begin{equation}\Psi_{l_0}([0])=\Psi_{l_0}([a_0])\not=\Psi_{l_0}([b_0])\end{equation}
holds, where $\Psi_{l_0}$ denotes the $l_0$-th component of $\Psi$. 
Then, because $\Psi$ is a $C(T;\mathbb R)$-injection 
from $(C(T;\mathbb R))/(\mathbb Z a_0)$ to $(C(T;\mathbb R))^k$, as we put 
$$\Phi (x) \ := \ \Psi_{l_0}([x]) \ \ \ \ \ \ (x \in C(T;\mathbb R)),$$ 
the map  
$$x \, \in \, C(T;\mathbb R) \ \ \ \mapsto \ \ \ \Phi (x) \, \in \, C(T;\mathbb R)$$
is a $C(T;\mathbb R)$-map and, especially, 
Frechet derivatives of $\Phi$ are $C(T;\mathbb R)$-linear transformations 
on $C(T;\mathbb R)$. That is, 
\begin{equation}\Phi^\prime (x) \, \in \, C(T;\mathbb R) 
\ \ \ \ \ \ (x \in C(T;\mathbb R))\end{equation}
holds, because, in general, 
if $A$ is a commutative ring with the unity $1_A$ 
and $L$ is an $A$-linear transformation on $A$, then 
$$L(1_A) \ \in \ A$$
and 
$$L(x)=L(x\cdot 1_A)=x\cdot (L(1_A))=(L(1_A))\cdot x  \ \ \ \ \ \ (x \in A)$$
hold. 
Now, from (2), 
there exists $t_0\in T$ such that  
\begin{equation}(\Phi(0))(t_0)=(\Phi(a_0))(t_0)\not=(\Phi(b_0))(t_0)\end{equation}
holds. Then, from (1), 
Case $1: \, ``0=b_0(t_0)\mbox{''}$ 
or Case $2: \, ``a_0(t_0)=b_0(t_0)\mbox{''}$ holds. 

First, we consider Case $1$. Let $0=b_0(t_0)$. 
Let 
$$F_1(s) \ := \ s b_0 \ \ \ \ \ \ (s\in [0,1]).$$
Then, the map 
$$s \, \in \, [0,1] \ \ \ \mapsto \ \ \ F_1(s) \, \in \, C(T;\mathbb R)$$ 
is a line from $0$ to $b_0$ in $C(T;\mathbb R)$ and 
$$(F_1)^\prime(s) \, = \, b_0 \ \ \ \ \ \ (s\in [0,1])$$
holds. Hence, 
$$\Phi(b_0)-\Phi(0)$$
$$= \ \Phi(F_1(1))-\Phi(F_1(0))$$
$$= \ \int_0^1 \, (\Phi^\prime(F_1(s))) \, ((F_1)^\prime(s)) \, ds $$
$$= \ \int_0^1 \, (\Phi^\prime(F_1(s))) \, b_0 \, ds$$
holds. So, in virtue of (3),  for any $t\in T$, 
%%%%%%%%%%%%%%%%%%%%%%
%%%%%%%%%%%%%%%%%%%%%%%%
%%%%%%%%%%%%%%%%%%%%%%%%%%
\newpage
%%%%%%%%%%%%%%%%%%%%%%%%%%%%5
%%%%%%%%%%%%%%%%%%%%%%%%%%%%
%%%%%%%%%%%%%%%%%%%%%%%%%%%%%%
$$(\Phi(b_0))(t)-(\Phi(0))(t)$$
$$= \ \int_0^1 \, ( \, (\Phi^\prime(F_1(s))) \, b_0 \, ) \, (t) \, ds$$
$$= \ \int_0^1 \, ((\Phi^\prime(F_1(s)))(t)) \cdot (b_0(t)) \, ds$$
holds. However, because of $0=b_0(t_0)$, 
$$(\Phi(b_0))(t_0)-(\Phi(0))(t_0)=0$$ holds. It contradicts (4).  

Next, we consider Case $2$. Let $a_0(t_0)=b_0(t_0)$. 
Let 
$$F_2(s) \ := \ a_0 \, + \, s (b_0-a_0) \ \ \ \ \ \ (s\in [0,1]).$$
Then, similarly, because 
$$\Phi(b_0)-\Phi(a_0)$$
$$= \ \int_0^1 \, (\Phi^\prime(F_2(s))) \, (b_0-a_0) \, ds$$
holds, in virtue of (3) and $a_0(t_0)=b_0(t_0)$, 
$$(\Phi(b_0))(t_0)-(\Phi(a_0))(t_0)$$
$$= \ \int_0^1 \, ((\Phi^\prime(F_2(s)))(t_0)) \cdot (b_0(t_0)-a_0(t_0)) \, ds$$
$$=0$$ holds. It contradicts (4).  
\hfill 
$\blacksquare$

\vspace*{0.8em}

\noindent 
{\sf Comment} : \ 
Kasuya suggested that an $\mathbb R^n$-manifold 
and a $\mathbb C^n$-manifold might be related 
to a differential web and a holomorphic web, respectively.  
He proposed a candidate for a compact $\mathbb C^2$-manifold $N$ 
such that for any $\mathbb C$-manifolds 
$M_1$ and $M_2$, $N$ can not be embedded 
in $M_1\times M_2$ as a $\mathbb C^2$-submanifold. 
The construction of our example was inspired by this proposal. 
\mbox{} 
\hfill ---

\vspace*{0.8em} 
\vfill
%\vspace*{3.2em}

\noindent  
Acknowledgment: \ \ \ 
As in Comment, Professor Naohiko Kasuya suggested it. 
This work was supported by JSPS KAKENHI Grant Number JP16K05245.

%\newpage

%\noindent 
%{\bf Oral} : \ 

%If there is not much previous research, it will be either

%\noindent
%Case 1: It has not been considered seriously.

%\noindent
%or

%\noindent
%Case 2: It is difficult to find related interesting problems. Or, even if you find a problem, it is difficult to get even a ``partial result''.

%By the way, manifolds on a (real or complex) commutative Banach algebra are very natural objects. Nevertheless, even an ordinary mathematician like me can get some ``slight partial results'' with just a few thoughts, and is likely to hit the case 1.

%Now, it is great to create a novel breakthrough for a subject that has been studied extensively, but this is difficult. If you are a mediocre mathematician like me, 
%you will likely be able to do a much better study on a subject that has not yet been considered.

%But, there is one problem here. Because the research field has not been established, even if a result is obtained, there are not many people who immediately evaluate it. Therefore, it is not recommended that people who want to acquire academic posts in the future do their main research. To put it the other way around, it is highly recommended for permanent employees.
%\mbox{} \hfill ---

%\vspace*{0.8em}

%\vfill
%\vspace*{3.2em}

%\noindent  
%Acknowledgment: 

%As in Comment, Professor Naohiko Kasuya suggested it. 

%This work was supported by JSPS KAKENHI Grant Number JP16K05245. 
%$ \mbox{ \ \ \ \ \ \ } $ 
%\hfill ---

\newpage 

{\bf References}

%\vspace*{0.8em}

%[1] H. Cartan and P. Thullen, Zur Theorie der Singularitaten der Funktionen mehrerer %komplexen Veranderlichen (German), {\it Math. Ann.}, 106 (1932),
%617-647.

[1] M. Abel and M. Abel, 
Modified definition of an $A$-bundle and a version 
of the Serre-Swan-Mallios theorem 
for general topological algebras, 
{\it Rend. Circ. Mat. Palermo}, 58 (2009), 345-360. 

[2] B. W. Glickfeld, The Riemann sphere of a commutative Banach algebra, 
{\it Trans. Amer. Math. Soc.}, 134 (1968), 1-28.

%[2] K. R. Goodearl, Cancellation of low-rank vector bundles, 
%{\it Pacific J. Math.}, 113 (1984), 289-302.

%[3] L. Hormander, 
%$L^2$ estimates and existence theorems for the $\overline \partial$ operator,
%{\it Acta Math.}, 113 (1965), 89-152.

[3] S. Kobayashi, Manifolds over function algebras and mapping spaces,
{\it Tohoku Math. J.}, 41 (1989), 263-282.

%[5] L. Lempert, The Dolbeault complex in infinite dimensions, 
%{\it J. Amer. Math. Soc.}, 11 (1998), 485-520.

%[4] E. R. Lorch, The theory of analytic functions in normed Abelian vector
%rings, {\it Trans. Amer. Math. Soc.}, 54 (1943), 414-425.

[4] A. Mallios and E. E. Rosinger, Space-time foam dense singularities
and de Rham cohomology, {\it Acta Appl. Math.}, 67 (2001), 59-89.

%[6] P. Manoharan, A nonlinear version of Swan's theorem, 
%{\it Math. Z.}, 209
%(1992), 467-479.

%[7] P. Manoharan, Generalized Swan's theorem and its application, 
%{\it Proc. Amer. Math. Soc.}, 123 (1995), 3219-3223.

[5] P. Manoharan, A characterization for spaces of sections, 
{\it Proc. Amer. Math. Soc.}, 126 (1998), 1205-1210.

%[9] A. S. Morye, Note on the Serre-Swan theorem, 
%{\it Math. Nachr.}, 286
%(2013), 272-278.

%[12] T. Ohsawa and K. Takegoshi, On the extension 
%of $L^2$ holomorphic
%functions, {\it Math. Z.}, 195 (1987), 197-204.

[6] M. H. Papatriantafillou, Partitions of unity on $A$-manifolds, 
{\it Internat. J. Math.}, 9 (1998), 877-883.

%[11] R. G. Swan, Vector bundles and projective modules, {\it Trans. Amer.
%Math. Soc.}, 105 (1962), 264-277.

%[12] L. N. Vaserstein, Vector bundles and projective modules, 
%{\it Trans. Amer. Math. Soc.}, 294 (1986), 749-755.

[7] H. Yagisita, Finite-dimensional complex manifolds on commutative Banach algebras and continuous families of compact complex manifolds, 
{\it Complex Manifolds}, 6 (2019), 228-264. 

[8] H. Yagisita, Cartan-Thullen theorem for a $\mathbb C^n$-holomorphic convexity  
and a related problem, {\it preprint}.

[9] H. Yagisita, A remark on locally direct product subsets 
in a topological Cartesian space, {\it preprint}.

[10] H. Yagisita, A manifold on the real commutative Banach algebra $C([0,1];\mathbb R)$
that cannot be embedded in the finite-dimensional Euclidean space
$C([0,1];\mathbb R^n)$, {\it preprint}.

\vfill 
\vspace*{0.8em} 
\noindent 
The related literature: \ 
``https://www.researchgate.net/profile/Hiroki\_Yagisita''

%The revised version 
\noindent 
%may have been put in 
%``https://www.researchgate.net/profile/Hiroki\_Yagisita''. 

%$ \mbox{ \ \ \ \ \ \ } $ 
%\hfill ---

%\noindent Running title: Expanding fronts in an anisotropic diffusion equation  
%\vspace*{0.8em}

\end{document}